\documentclass[conference]{IEEEtran}
 \usepackage[dvips]{graphicx}
 \DeclareGraphicsExtensions{.eps}
\begin{document}
%
\title{Free Final-Time Optimal Control for HIV Viral Dynamics}

\author{\IEEEauthorblockN{Gaurav Pachpute}
\IEEEauthorblockA{Department of Mathematics,\\
Indian Institute of Technology Guwahati\\
Guwahati, Assam 781039\\
Email: g.pachpute@iitg.ernet.in}
\and
\IEEEauthorblockN{Kriti Saxena}
\IEEEauthorblockA{Electronic and Instrumentation Engineering,\\
BITS-Pliani, Dubai Campus,\\
Dubai, UAE\\
Email: kriti217@gmail.com}
}


%


\maketitle

\begin{abstract}
In this paper, we examine a well-established model for HIV wild-type infection. The algorithm for steepest descent method for fixed final-time is stated and a modified method for free final-time is presented. The first type of cost functional considered, seeks to minimize the total time of therapy. An easy implementation for this problem suggests that it can be effective in the early stages of treatment as well as for individual-based studies, due to the ``hit first and hit hard'' nature of optimal control.

An LQR based cost functional is also presented and the solution is found using steepest descent method. It suggests that the optimal therapy must remain high until the patient shows signs of recovery after which, the therapy gradually decreases. This is in line with the biomedical philosophy. Solution to a modified problem which includes a weight for total time is approximated using the modified algorithm. It shows a considerable drop in the total period. We conclude that, a decreased and optimized therapy period can help us increase efficiency as well as the turnover rate for patient care.
\end{abstract}


%
\IEEEpeerreviewmaketitle
\def\Jc{\mathcal{J}}
\def\Hc{\mathcal{H}}
\def\Lc{\mathcal{L}}
\def\Kc{\mathcal{K}}
\def\xbf{\mathbf{x}}
\def\ubf{\mathbf{u}}
\def\fbf{\mathbf{f}}
\def\lambf{\mathbf{\lambda}}
\def\Jbf{\mathbf{J}}
\def\SS{\scriptsize}
\def\TM{\!\!$^{\mbox{\texttt{\SS TM}}}$~}

\section{Introduction}
The World health organization (WHO) classifies human immunodeficiency virus (HIV) as a pandemic. HIV progresses into acquired immunodeficiency syndrom (AIDS) which leads to  immune system failure. By 2006, AIDS has claimed over 25 million lives and over 0.6\% of the world population is infected with HIV\cite{WHO06}.

HIV is responsible for the selective depletion of $CD4^{+}$ cells, also known as helper T cells. As they are essential to immune regulation, such depletion leads to an adverse effect on the immune system functioning \cite{perelson93}. Even though HIV is rarely fatal by itself, it increases vulnerability towards infections and malignancies \cite{stengel08}.
There is no known cure for HIV/AIDS as the available drug regimens fail to eliminate HIV strains in overall population \cite{dimascio04}. The antiretroviral treatments available are mostly administered in ``drug-cocktail'' form, also known as highly active antiretroviral therapy (HAART). Even though HAART has been highly effective in treating HIV/AIDS, it leads to numerous side-effects from hepatitis, liver failure, cardiovascular malfunction, pancreas damage to nausea, diarrhea and depression \cite{stengel08}. In most patients, the therapy is long-term and therefore, the search for optimal therapy stems from the notion of maintaining a balance between the disease and drug side-effects to minimize patient suffering.

Several models have been proposed for different types of HIV therapies. Many of these models study the host-pathogen reaction for the virus strain HIV-1 through medical studies or numerical simulations \cite{stengel08}. McLean and Nowak \cite{mclean92} address the appearance of AZT-resistant strains and its effects on treatment. Agur \cite{agur89} and Agur and Cojucaru \cite{agur93} examine the effects of chemotherapy on uninfected cells. Much of the research is also dedicated to acquired infections in an HIV patient \cite{stengel08}. Compartment infections are studied and modeled in \cite{perelson97} and \cite{bajaria02}. Microphage infections are addressed in \cite{wodarz99}, whereas, transient viremia is discussed in \cite{jones05}.

In this paper, we discuss the wild-type infection model given in \cite{stengel08}. This model is a combination of several of the previous models, including, \cite{kirschner97}, \cite{perelson89}, and \cite{perelson93}. This model encapsulates four different therapies : protease inhibitors, HIV fusion inhibitors, T cell enhancer and HIV reverse transcriptase inhibitor. In this paper, we study the treatment with protease inhibitors.

\section{HIV Model for Wild-Type Infection}
This paper uses the model given by Stengel \cite{stengel08} for wild-type infectious HIV. The model is given by four coupled-ODEs and their elements are the HIV particles ($x_{1}$), uninfected $T_{h}$ cells ($x_{2}$), proviral $T_{h}$ cells ($x_{3}$) and productively infected  $T_{h}$ cells ($x_{4}$) (all in per mm$^{3}$). The different therapies are protease inhibitors $(u_{1})$, HIV fusion inhibitors $(u_{2})$, T cell enhancer $(u_{3})$ and HIV reverse transcriptase inhibitor $(u_{4})$. The model is given by the following ODEs,

\begin{eqnarray}
\dot{x_{1}} &=& -a_{1} x_{1}-a_{2}x_{1}x_{2}(1-u_{2})+a_{3}a_{4}x_{4}(1-u_{1}) \nonumber\\
\dot{x_{2}} &=& \frac{a_{5}}{1+x_{1}}-a_{2}x_{1}x_{2}(1-u_{2})(1-u_{4})-a_{6}x_{2} \nonumber\\
& & +a_{7}\left(1-\frac{x_{2}+x_{3}+x_{4}}{a_{8}}\right)x_{2}(1+u_{3}) \nonumber\\
\dot{x_{3}} &=& a_{2}x_{1}x_{2}(1-u_{2})(1-u_{4})-a_{9}x_{3}-a_{6}x_{3} \nonumber\\
\dot{x_{4}} &=& a_{9}x_{3}-a_{4}x_{4}
\end{eqnarray}

The virus particles ($x_{1}$) have a death rate of $a_{1}$. $x_{1}$ infect the uninfected T cells ($x_{2}$) at rate $a_{2}$. This is partially blocked by HIV fusion inhibitors with efficacy $u_{2}$. A portion (1-$u_{1}$) of the production of $x_{1}$ from the productively infected cells ($x_{4}$) is blocked by the protease inhibitors with efficacy $u_{1}$. Birth rate of $x_{2}$ is proportional to $a_{5}$ and depends on the amount of virus present in the system. The conversion of $x_{2}$ to proviral $T_{h}$ cells $(x_{3})$ occurs at a rate $a_{2}$ and is targeted by HIV reverse transcriptase inhibitor $(u_{4})$. $x_{2}$ has a natural death rate of $a_{6}$ and a proliferation rate of $a_{7}$. Proliferation is enhanced by T cell enhancer $(u_{3})$. Proviral $T_{h}$ cells $(x_{3})$ are produced through infection. $x_{3}$ convert into $x_{4}$ at a rate $a_{9}$ and have a death rate of $a_{6}$. The infected cells $(x_{4})$ convert into virus $(x_{1})$ at rate $a_{4}$ with the ratio 1:$a_{3}$, completing the cycle.
The parameter values (from \cite{stengel08} and \cite{kirschner97}) are $a_{1}=2.4$, $a_{2}=2.4\times 10^{-5}$, $a_{3}=1200$, $a_{4}=0.24$, $a_{5}=10$, $a_{6}=0.02$, $a_{7}=0.03$, $a_{8}=1500$, $a_{9}=0.003$.

We define the vectors $\xbf=[x_{1}~x_{2}~x_{3}~x_{4}]^{\top}$ and $\ubf=[u_{1}~u_{2}~u_{3}~u_{4}]^{\top}$, then equation (1) can be written as $\dot{\xbf}=f(\xbf(t),\ubf(t),t)$. The initial conditions are given by $\xbf(0)=\xbf_{0}$. We assume that, the drug efficacy for drug $i$, $u_{i}$ varies between $[\underline{u_{i}}, \overline{u_{i}}]$ and as the efficacies are normalized this is a subset of $[0,1]$. In this paper, we discuss the therapy with protease inhibitors ($u_{1}$) and take $u_{2}=u_{3}=u_{4}=0$.

\section{Optimal Control Theory}
A general optimal control problem can be stated as a combination of a system of ordinary differantial equations (ODEs), which defines the dynamics of the state variables and controls along with a cost functional, which is to be minimized. The problem is defined over a time period $[t_{o},t_{f}]$, where $t_{f}$ can be predefined or free.  Consider a system given by $\dot{\xbf}=f(\xbf,\ubf,t)$, $\xbf(0)=\xbf_{0}$ and a cost functional given by,
\[\Jc=\Kc(\xbf(t_{f}),t_{f})+\int _{t_{o}}^{t_{f}} \Lc(\xbf(t),\ubf(t),t) dt. \]
If the final time is free, another constraint can be put on the final state of the system : $\Psi(\xbf(t_{f}),t_{f})=0$, where $\Psi$ is a smooth function.
\textit{Hamiltonian} for this system is defined by,
$\Hc(\xbf,\ubf,\lambda,t)=\Lc(\xbf,\ubf,t)+\lambda(t) f(\xbf,\ubf,t)$. The necessary conditions for a minimum are achieved by defining the modified cost functional $\tilde{\Jc}$ using lagrangian multipliers $\lambda(t)$ and $p$ as,
\begin{eqnarray}
\tilde{\Jc}&=&\Kc(\xbf(t_{f}),t_{f})+p \Psi(\xbf(t_{f}),t_{f}) + \nonumber\\
& &\int _{t_{o}}^{t_{f}} \left[\Lc(\xbf(t),\ubf(t),t) -\lambda^{\top}(t) (f(\xbf(t),\ubf(t),t)-\dot{x})\right] dt.\nonumber
\end{eqnarray}

The cost functional reaches minimum when the variation, $\delta \tilde{\Jc} =0$ with respect to all variables, $\lambda, p, \xbf, \ubf$ and $\xbf(t_{f})$ (also, $t_{f}$, if final-time is free). Using this argument, if $\xbf^{*}(t)$ and $\ubf^{*}(t)$ are the optimal state and control vectors for time $t$ respectively then we can state the Pontryagin's necessary conditions as following,
\begin{enumerate}
\item $\Psi(\xbf^{*}(t_{f}),t_{f})=0$ (Final State Constraint)
\item $\dot{\xbf}^{*}=\nabla_{\lambda} \Hc=f(\xbf^{*}(t),\ubf^{*}(t),t) $ (State Equation)
\item $\dot{\lambda}=-\nabla_{\xbf^{*}} \Hc$ (Costate Equation)
\item $\nabla_{\ubf^{*}} \Hc =0 $ (Optimal Control)
\item $\xbf(t_{o})=\xbf_{0}$\\
$\lambda(t_{f})=\nabla_{x(t_{f})} (\Kc+p^{\top} \Psi)$\\
$\Hc(t_{f})=-\nabla_{t_{f}} (\Kc+p^{\top} \Psi)$ (Boundary Conditions)
\end{enumerate}
As for the sufficient conditions, where one is concerned with the local minimality, it is given by $\nabla_{u^{*}}^{2} \Hc > 0$ (Positive definite). For global minimality, one must consider all possible controls over time $[t_{o},t_{f}]$. A Hamiltonian  convex in $\ubf$ can ensure the global minimality.

If the final time is free, only the first two boundary conditions are valid. Furthermore, if there is no final state constraint, the second boundary condition becomes, $\lambda(t_{f})=\nabla_{x(t_{f})}\Kc$.

One of the most widely used cost functional is the Linear-Quadratic Regulator (LQR) given by,
\[\Jc=\frac{1}{2}\xbf^{\top}(t_{f})S\xbf(t_{f}) + \int_{t_{o}}^{t_{f}} \frac{1}{2}\left[\xbf^{\top}(t)Q\xbf(t) + \ubf^{\top}(t)R\ubf(t)\right] dt.\]
The matrices, $S, Q, R$ are the weights for the final state configuration, integral value of state and control, respectively. If the final time is free, one may add a scalar $T$ in the integral to include the weight for total time.

\subsection{Algorithm for fixed final-time optimal control problems}
As it can be seen from the necessary conditions, the optimal control system comprises a forward state equation and a backward costate equation. For nonlinear models the algorithms are essentially iterative. One of the very effective of these algorithms is the steepest descent method. The solutions in this paper are approximated using this algorithm \cite{pachpute11}.

\begin{enumerate}
\item Approximation starts with an initial guess $\ubf(t)=\ubf^{(0)}(t)$.
\item $\xbf^{(0)}(t)$ is integrated numerically over time $[t_{o},t_{f}]$ using control $u^{(0)}(t)$ with $\xbf(0)=\xbf_{0}$ as the boundary condition.
\item Similarly, the costate equation $\dot{\lambda}=-\nabla_{\xbf^{(0)}} \Hc$ is numerically solved with $\lambda(t_{f})=\nabla_{x(t_{f})} \Kc$ as the boundary condition.
\item The gradient $\nabla_{\ubf^{(0)}} \Hc$ is evaluated using numerical techniques and the updated control for the next iteration is given by,
    \[\ubf^{(1)}(t)=\ubf^{(0)}(t)-\tau(0)\nabla_{\ubf^{(0)}} \Hc.\]
    Similarly, the k+1-th iteration in terms of $u^{k}(t)$ is given by,
    \[\ubf^{(k+1)}(t)=\ubf^{(k)}(t)-\tau(k)\nabla_{\ubf^{(k)}} \Hc.\]
\item This iterative procedure is repeated until some convergence criterion such as, $||\ubf^{(k+1)}-\ubf^{(k)}||<tol$ or $\Jc^{(k)}-\Jc^{(k+1)}<tol$ is satisfied.
\end{enumerate}

\subsection{Modified algorithm for free-final time problems}
If the final time is free, the algorithm can be modified to encapsulate it by updating the control as well as the grid size for numerical integration in every iteration. For a small variation $\Delta t_{f}$, one can write $\Delta \Jc=\frac{\delta \Jc}{\delta t_{f}} \Delta t_{f}$. For LQR, this can be written as,

\begin{eqnarray}
\label{dotlqr}
\Delta \Jc &=& \left[\xbf{\top}(t_{f})S \nabla_{t_{f}} \xbf(t_{f})\right]\Delta t_{f} + \nonumber \\
& &\frac{1}{2}\left[\xbf^{\top}(t_{f})Q\xbf(t_{f}) + \ubf^{\top}(t_{f})R\ubf(t_{f})+T\right]\Delta t_{f} \nonumber\\
&=&\xbf{\top}(t_{f})Sf(x(t_{f}),u(t_{f}),t_{f})\Delta t_{f} + \nonumber\\
& & \frac{1}{2}\left[\xbf^{\top}(t_{f})Q\xbf(t_{f}) + \ubf^{\top}(t_{f})R\ubf(t_{f})+T\right]\Delta t_{f} \nonumber\\
& &
\end{eqnarray}

Suppose, one expects to reduce $\Jc$ by 10\%, then one can change $t_{f}$ using the following equation,
\[\Delta t_{f} = t^{(k+1)}_{f}-t^{(k)}_{f}= -\frac{0.1}{\frac{\delta \Jc}{\delta t_{f}}},\]
where, $t^{(k)}_{f}$ represents the k-th iteration. Note that, if the number of grid points is $n$ and the grid length at the k-th iteration is $h^{(k)}$ then $t_{f}^{(k)}=(n-1)h^{(k)}$. Therefore, at every iteration the grid length can be updated as,
\[h^{(k+1)}=h^{(k)}-\frac{\tau_{h}(k)}{(n-1) \frac{\delta \Jc}{\delta t_{f}}}.\]

This condition corresponds to the third boundary condition specified in the last section, i.e., $\frac{\delta \Jc}{\delta t_{f}}=\Hc(t_{f})+\nabla_{t_{f}} \Kc$. However, this algorithm cannot be effective when the problem is specified with a final state constraint ($\Phi$).

\section{Results and Discussions}
\subsection{Minimum-time problem}
In minimum-time problems, the cost functional is defined as $\Jc=t_{f}-t_{o}=\int_{t_{o}}^{t_{f}} dt$. They are often presented with an acceptable set $\Omega$ of final state variables, i.e., the value of $\Phi$ over this set is zero. The problem specified in the next paragraph can be solved using the iterative method described in the last section by continuing the process until the system either enters the set $\Omega$ or leaves it, but as there is no incentive for the control to be any less than its maximum value, another way to define these problems is to find the minimum time it takes for the system to enter set $\Omega$. Using this notion, we can solve the problem with an error of the order of the grid size used for numerical integration.

\begin{figure}[!t]
\centering
\includegraphics[width=3.5in]{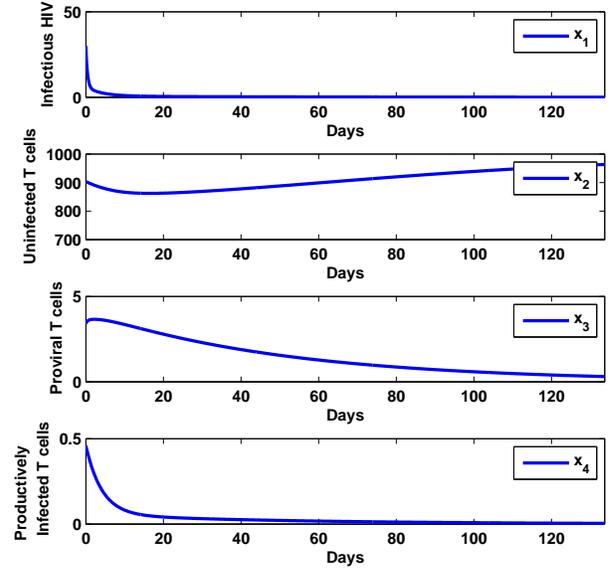}
\caption{States for min-time problem}
\label{fig1}
\end{figure}

\begin{figure}[!t]
\centering
\includegraphics[width=3.5in]{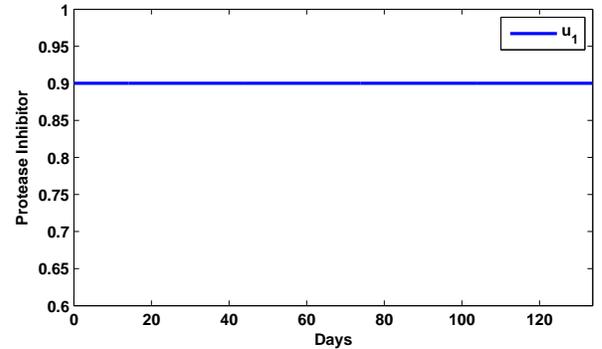}
\caption{Control for min-time problem}
\label{fig2}
\end{figure}

For the system at hand, the set $\Omega$ is taken as the values of virus particles below a certain level, often called the detection level. We take the detection level to be 0.05 \cite{barletta04}. The minimum-time problem can be stated as :
To find an admissible control $u_{1}^{*} : [0,t_{f}]\rightarrow[\underline{u_{1}}, \overline{u_{1}}]$ and the final time $t_{f}$, which for the system $\dot{\xbf}=f(\xbf(t),\ubf(t),t), \xbf(0)=\xbf_{0}$ minimizes the cost functional $\Jc=t_{f}-t_{o}=t_{f}$ such that the final virion population is below the detection level. Here, the boundary conditions are taken as $x_{1}(0)=30.0, x_{2}(0)=904, x_{3}(0)=3.4, x_{4}(0)=0.46$ and $\overline{u_{1}}=0.9$, \cite{brogan10} for pharmacokinetic models for HIV suggest that the efficacy decays exponentially and a perfect efficacy throughout the treatment period is unlikely to be achieved. The solution to this problem is presented in Figure \ref{fig1} and Figure \ref{fig2}. The results suggest that it takes about 19 weeks for the virus to climb down below the detection level.

This type of method can be used for finding the optimal therapy for curable or highly manageable infections or for drugs with minimal side-effects, where a full-blast treatment is not detrimental, since the main cause of suffering originates from the time the infection is active in the system. This method can also be useful in an early detection or treatment for the disease, as we shall see in the next subsection that the therapy remains at its highest under these conditions. The maximum efficacy can be averaged over one period of administering drugs for a full-blast therapy. The easy implementation of the aforementioned procedure leaves room for studies based on individual patient scenarios \cite{pachpute11(2)}.

\subsection{Fixed final-time problem}
In this section, as a prerequisite to the free final-time problem, we develop a fixed final-time cost functional (after \cite{stengel08}) which sheds light on some of the key features of the optimal control. We take $x_{1}(0)=4.9, x_{2}(0)=904, x_{3}(0)=0.34, x_{4}(0)=0.42$ and $t_{f}=500$ days. The cost functional is given by,
\[\Jc=\frac{1}{2}\xbf^{\top}(t_{f})S\xbf(t_{f}) + \int_{t_{o}}^{t_{f}} \frac{1}{2}\left[\xbf^{\top}(t)Q\xbf(t) + \ubf^{\top}(t)R\ubf(t)\right] dt,\]
where, S, Q and R are diagonal matrices. We wish to design a problem that minimized the terminal and integral values of proviral T cells ($x_{3}$), productively infected cells ($x_{4}$) and the virus ($x_{1}$) as well as the integral value of the control ($u_{1}$).

\begin{figure}[!t]
\centering
\includegraphics[width=3.5in]{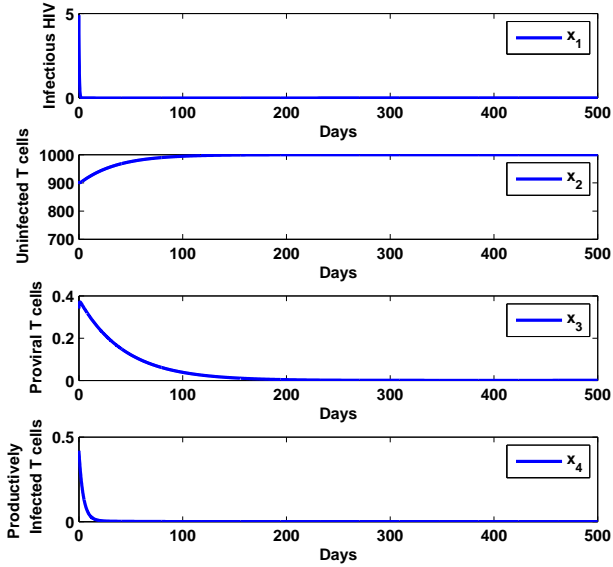}
\caption{States for fixed final-time problem}
\label{fig3}
\end{figure}

\begin{figure}[!t]
\centering
\includegraphics[width=3.5in]{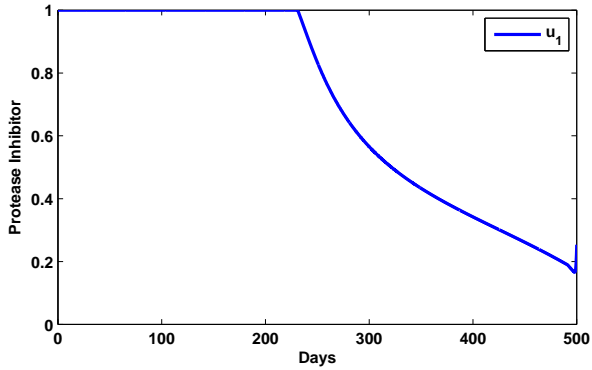}
\caption{Control for fixed final-time problem}
\label{fig4}
\end{figure}

Keeping these goals in mind, we take the weights to be, $S_{11}=S_{33}=S_{44}=10^{3}$, $Q_{11}=Q_{33}=Q_{44}=10^{3}$ and $R_{11}=0.01$ \cite{stengel08}. The solution is approximated by steepest decent method described in the previous section and is presented in Figure \ref{fig3} and Figure \ref{fig4}.

The results suggest that the optimal therapy ($u^{*}$) stays at its maximum for about 230 days and then gradually declines over the rest of the period. This is in sync with
the bio-medical philosophy of ``hit first and hit hard'' \cite{stengel08}, i.e., to treat the disease with high dosage in the early phase of the treatment and once a physical recovery takes place, gradually decrease the dosage. The physical recovery is represented by the uninfected $T_{h}$ cell count, which rises to almost its highest in the first 230 days \cite{pachpute11}. Similarly, the proviral count and the productively infected T cells have gone down substantially. Simulations run over a larger or smaller periods along with varied values of S,Q, and R suggest that the period of high efficacy does not vary by much as long as the values of R stay relatively low. This suggests that when treating with drugs with minimal side effects, the length of this period stays approximately the same \cite{pachpute11}.

\begin{figure}[!t]
\centering
\includegraphics[width=3.5in]{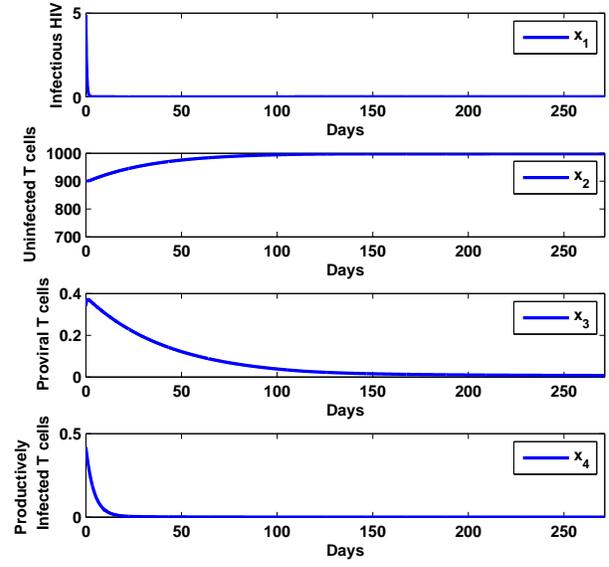}
\caption{States for free final-time problem}
\label{fig5}
\end{figure}

\begin{figure}[!t]
\centering
\includegraphics[width=3.5in]{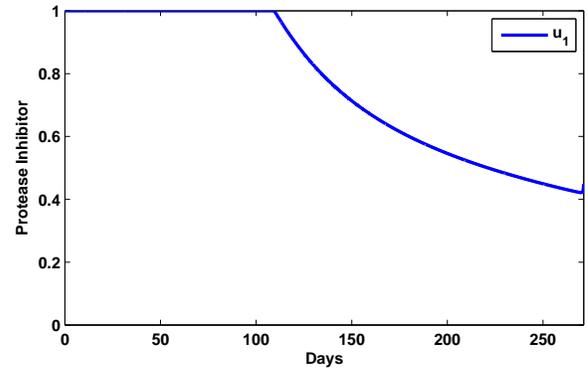}
\caption{Control for free final-time problem}
\label{fig6}
\end{figure}

This kind of studies of optimal therapy are very common and successful, as they reflect the well-accepted conceptions of the bio-medical field. They do, however, fail to address the time window for which the therapy continues. This is not reflected through the drug side-effects, especially, if they are considerably high, for the optimal efficacy changes with the period of treatment. Therefore, it is important to take time into account when trying to minimize the cost functional.

\subsection{Free final-time problem}
In this section, we present the model for a free final-time LQR. We state the optimal control problem as the following,

To find an optimal control $u_{1}^{*}:[0,t_{f}]\rightarrow[0,1]$ which, for the system $\dot{\xbf}=f(\xbf(t),\ubf(t),t)$, minimizes the cost functional,
\begin{eqnarray}
\Jc&=&\frac{1}{2}\xbf^{\top}(t_{f})S\xbf(t_{f}) +T t_{f} \nonumber \\
& &+\int_{0}^{t_{f}} \frac{1}{2}\left[\xbf^{\top}(t)Q\xbf(t) + \ubf^{\top}(t)R\ubf(t)\right] dt, \nonumber
\end{eqnarray}
where, the weights S, Q and R serve the same purpose as before and the added term $T t_{f}$ adds a cost for the duration of therapy. The values used are, $S_{11}=S_{33}=S_{44}=10^{3}$, $Q_{11}=Q_{33}=Q_{44}=10^{3}$, $R_{11}=0.01$ and $T=0.001$. The initial state of the system is taken to be the same as for the previous problem, i.e., $x_{1}(0)=4.9, x_{2}(0)=904, x_{3}(0)=0.34, x_{4}(0)=0.42$.

To solve this problem, the algorithm described in section 2 is implemented. $\frac{\delta \Jc}{\delta t_{f}}$ is given by equation \ref{dotlqr} and the grid length $h$ is updated after every few iterations. As mentioned in section 2, this method tries to satisfy one of the boundary conditions, namely, $\Hc(t_{f})+\nabla_{t_{f}} \Kc=0$. An implementation in MatLab$^{TM}$ shows a considerable drop in this value from the initial state. The results are presented in Figure \ref{fig5} and Figure \ref{fig6}.

As the results suggest, we notice a substantial drop in the treatment window from 500 days to about 275 days. However, the essential features of the therapy described in the last subsection are retained in Figure \ref{fig6}. The optimal therapy stays at its highest for the first 120 days and then shows a gradual decrease. This, as noted before, coincides with the restored T cell count and the decreased adverse element counts in the system. As given in the last subsection, this period of high efficacy varies (from 230 days to 120 days) when treated over a different time window. The state variables show similar progression over time under the two optimal therapies. The method decreases the total period of therapy while keeping the state variables as intact as possible.

\section*{Conclusion}
In this paper, we consider a well established model for wild-type HIV infection. A common LQR optimal therapy is approximated with the steepest descent method. We, also, present a modified algorithm for approximating the solution for free final-time problems.

A minimum-time problem is stated and the solution is approximated with a non-iterative method. We conclude that, due to the simple implementation and the ``hit first and hit hard'' nature of the optimal therapy, this kind of studies can be very useful in the first phase of therapy and also provide an advantage for individual-based optimal-therapy considerations.

The second cost functional is an LQR defined over a fixed time period, which indicates similarities of the optimal control with biomedical philosophies. The optimal therapy shows a period of high efficacy in the first phase and a gradual decrease once there is physical progress in patient. A modified problem with free final-time is approximated using the algorithm mentioned earlier. It shows a considerable drop in the total period. An optimized time window for treatment can help us treat patients more efficiently. More number of patients can be treated if the optimal period for therapy is lower.

The method presented in this paper is quite similar to methods already in use. The need for optimizing the therapy period is quite essential as a higher or lower time period suggests a non-optimal exposure to either the disease or therapy.

\end{document}